\documentclass{article}

\begin{document}

\title{b-stability and blow-ups}
\author{S. K. Donaldson}
\maketitle
\newcommand{\hTr}{\widehat{{\rm Tr}}}
\newcommand{\hDim}{\widehat{{\rm Dim}}}
\newcommand{\hVol}{\widehat{{\rm Vol}}}
\newcommand{\cU}{{\cal U}}
\newcommand{\hatX}{\hat{X}}
\newcommand{\osigma}{\overline{\sigma}}
\newcommand{\otau}{\overline{\tau}}
\newcommand{\os}{\overline{s}}
\newcommand{\ord}{{\rm ord}}
\newcommand{\cF}{{\cal F}}
\newcommand{\cZ}{{\cal Z}}
\newcommand{\hpi}{\hat{\pi}}
\newcommand{\cI}{{\cal I}}
\newcommand{\cE}{{\cal E}}
\newcommand{\hcX}{\hat{{\cal X}}}
\newcommand{\wW}{\widehat{W}}
\newcommand{\tg}{\tilde{g}}
\newcommand{\cA}{{\cal A}}
\newcommand{\bQ}{{\bf Q}}
\newcommand{\cX}{{\cal X}}
\newcommand{\bC}{\mbox{${\bf C}$}}
\newcommand{\bR}{\mbox{${\bf R}$}}
\newcommand{\Var}{{\rm Var}}
\newcommand{\Av}{{\rm Av}}
\newcommand{\Vol}{{\rm Vol}}
\newcommand{\Dim}{{\rm Dim}}
\newcommand{\cO}{{\cal O}}
\newcommand{\cW}{{\cal W}}
\newcommand{\cL}{{\cal L}}
\newcommand{\Tr}{{\rm Tr}}
\newcommand{\Zmax}{Z_{{\rm max}}}
\newcommand{\Zmin}{Z_{{\rm min}}}
\newcommand{\Ch}{{\rm Ch}}
\newcommand{\bP}{\mbox{${\bf P}$}}
\newcommand{\uA}{\mbox{${\underline{A}}$}}
\newcommand{\uM}{\mbox{${\underline{M}}$}}
\newcommand{\um}{\mbox{${\underline{m}}$}}
\newcommand{\ur}{\mbox{${\underline{r}}$}}
\newtheorem{Goal}{Goal}
\newtheorem{question}{Question}
\newtheorem{thm}{Theorem}
\newtheorem{prop}{Proposition}
\newtheorem{lem}{Lemma}
\newtheorem{defn}{Definition}
\newtheorem{cor}{Corollary}
\newcommand{\oK}{\overline{K}}
\newcommand{\dbd}{\sqrt{-1} \partial\overline{\partial}}
\newcommand{\ulambda}{\underline{\lambda}}
\newcommand{\olambda}{\overline{\lambda}}
\newcommand{\Riem}{{\rm Riem}}
\newcommand{\Ric}{{\rm Ric}}


     \ \ \ \ \ \ \ \ \ \ \ \ \ \ \ \ \ \ \ \ \ \ \ {\it Dedicated to Professor V. V. Shokurov}
     
     \

     \

\section{Introduction}
In \cite{kn:D2} the author introduced a notion of \lq\lq b-stability''.
This is a variant of more standard notions of stability, designed to get around certain difficulties in the proof of the well-known conjectures relating stability to the existence of K\"ahler-Einstein metrics (in the case of Fano manifolds). The purpose of this article is to make progress towards a proof that a manifold which admits a K\"ahler-Einstein metric is b-stable. The method is a variant of the one introduced by Stoppa \cite{kn:S1}, \cite{kn:S2}.

Let $X$ be a compact complex manifold and $L\rightarrow X$ be a positive line bundle. The definition of b-stability begins by considering a power $s$ such that sections of $L^{s}$ give a projective embedding  $X\subset \bP $ and degenerations
$\pi: \cX\rightarrow \Delta $ with $\cX\subset \bP\times \Delta$ such that $\pi^{-1}(t)\cong X$ for $t\neq 0$ but with $W=\pi^{-1}(0)$ not isomorphic to $X$. One then goes on to consider powers $sp$ and certain birational modifications of $\cX$. The main restriction we make in the present paper is to confine our  attention throughout to {\it equivariant} degenerations or \lq\lq test configurations''. Thus we have maps $g_{t}:\cX\rightarrow \cX$ for $t\neq 0$ with $\pi(g_{t}(x))= t \pi(x)$.  For our purposes we can also fix $s=1$. So we now consider a family $\cX\subset \bP\times \bC$ invariant under a $\bC^{*}$ action defined by the standard action on $\bC$, and a linear  action on $\bP$.  We make the following  assumptions
\begin{enumerate}
\item $H^{0}(X,L)$  generates the ring $\bigoplus_{k} H^{0}(X, L^{k})$.
\item $\pi:\cX\rightarrow \bC$ is a flat family, $\pi^{-1}(1)\cong X$ and $X\subset \bP$ is the embedding defined by the complete linear system $\vert L\vert$. 
\item $\pi^{-1}(0)=W$ is reduced and contains a preferred component $B$. We write
$W=B\cup R$ where $R$ is a union of one or more other components and set $D=B\cap R$.  
\item $B$ does not lie in any proper linear subspace in $\bP$.
\item For each $k$ the power $\cI_{B}^{k}$ of the ideal sheaf of $B$ in $\cX$ coincides with the sheaf of functions which vanish to order $k$ on $B$, in the sense of the valuation defined by $B$. 
\end{enumerate}

 Here the conditions (1)(2) are rather standard. Condition (4) will always hold in the applications such as described in \cite{kn:D2}. The assumption in condition (3) that $W$ is reduced can probably be relaxed. The last condition (5) is of a more technical nature. It is used to give a simple proof of the Lemma 4 below but requires clarification, see the discussion in Section 3. Notice that this condition (6) will hold if $B$ is a Cartier divisor in $\cX$, but that assumption would rule out many cases of interest (see the examples in Section 3 below). 

We now  recall the basic construction of \cite{kn:D2} in this situation.  Following a suggestion of Richard Thomas, we can present this construction  in a  slightly different (and probably more familiar) way than in \cite{kn:D2}. Let $\cL\rightarrow \cX$ be the pull-back of $\cO(1)$ and for $\mu>0$ let  $\Lambda^{(\mu)}$ be the sheaf of meromorphic functions with at worst poles of order  $\mu$ on $R$. So if $R$ is a Cartier divisor $\Lambda^{(\mu)}$ is just the sheaf of sections of the line bundle defined by $R$, raised to the power $\mu$.  For each positive integer $p$ we consider the inclusions
$$   i_{p,\mu}: H^{0}(\cX; \cL^{p}\otimes \Lambda^{(\mu-1)})\rightarrow H^{0}(\cX;\cL^{p}\otimes \Lambda^{(\mu)}) $$
and define $m(p)$ to be the largest value of $\mu$ such that $i_{p,\mu}$ is not an isomorphism. Now consider the  sheaf  $\cL'_{p}= \cL^{p}\otimes \Lambda^{(m(p))}$ over $\cX$. The sections of this define  a rational map from $\cX$ to $\bP_{p}\times \bC$, for a projective space $\bP_{p}=\bP(U_{p})$ and the image of this gives another flat family $\cX'_{p}$. This is essentially a restatement of the construction described in \cite{kn:D2}, as we will explain in Section 2 below. The upshot is that we get a new  degeneration $\pi^{(p)}:\cX'_{p}\rightarrow \bC$, embedded in $\bP_{p}\times \bC$ which is equivariant with respect to a $\bC^{*}$-action defined by a generator $A'_{p}\in {\rm End}(U_{p})$. In particular this gives a $\bC^{*}$-action on the central fibre $W'_{p}\subset \cX'_{p}$. (Note that strictly we should allow $\cX'_{p}$ and $ W'_{p}$ to be schemes, but this aspect will not really enter the discussion.)

\

We now come to another central topic in this paper: the {\it Chow invariant} of a projective variety with $\bC^{*}$-action. In general suppose we have an $n$-dimensional variety $V\subset \bP(\bC^{N})$ preserved by a $\bC^{*}$-action with generator $\alpha$. Let $\alpha_{k}$ be the generator of the induced action on $H^{0}(V;\cO(k))$. Then, for large $k$, the trace $\Tr(\alpha_{k})$ is given by a Hilbert polynomial of degree $n+1$ and we define $I(V)$ to be the leading term
\begin{equation}    I(V)= \lim_{k\rightarrow \infty} k^{-(n+1)} \Tr(\alpha_{k}). \end{equation} 
 Write $\Vol$ for the degree (or volume) of $V$. We define the Chow invariant of $V$ to be
\begin{equation} \Ch(V)= \frac{\Tr(\alpha)}{N}- \frac{I(V)}{\Vol}. \end{equation}
(Of course this depends on the given action, although we omit it from the notation.) We recall that a variety $Z\subset \bP^{N-1}$ is defined to be  {\it Chow
stable} if for all non-trivial flat equivariant degenerations of $Z$ contained in 
$\bC\times \bP^{N-1}$ the Chow invariant of the central fibre is strictly positive. 

\
We introduce some new definitions.
Let $q$ be a point in $X$ and $\hat{X}_{q}$ be the blow up of $X$ at $q$ with exceptional divisor $E\subset \hat{X}_{q}$ and corresponding line bundle $L_{E}$.
Let $\gamma, r$ be positive integers. We say that $(X,L)$ is {\it $(\gamma,r)$-stable} if for all points $q\in X$ the line bundle $L^{r\gamma} \otimes L_{E}^{-r}$
on the blow-up is very ample and if the corresponding  projective embedding of $\hat{X}_{q}$ is Chow stable. We say that $(X,L)$ is {\it bi-asymptotically stable} if there is a $\gamma_{0}$ and for each $\gamma\geq \gamma_{0}$ an
$r_{0}(\gamma)$ such that $(X,L)$ is $(\gamma,r)$-stable for $r\geq r_{0}(\gamma)$. (In \cite{kn:D2} we discussed a very similar notion of
\lq\lq $\overline{K}$-stability'', but on reflection the definition above seems to be the more relevant one.)

 The point of this definition is that we have
 \begin{thm}
 If the automorphism group of $(X,L)$ is finite and if $X$ admits a constant scalar curvature K\"ahler metric in the class $c_{1}(L)$ then $(X,L)$ is
bi-asymptotically stable.
\end{thm}

This is a standard result now, but we review the proof. 

\begin{itemize}
\item According to Arrezo and Pacard \cite{kn:AP}, if $\gamma$ is sufficiently large then for all $r>0$  there is a constant scalar curvature metric on the blow-up in the class $c_{1}(L^{r\gamma}\otimes L_{E}^{-r})$. 
\item Fix $\gamma$ and make $r$ large enough that $L^{r\gamma}\otimes L_{E}^{-r}$ is very ample. By the main result of \cite{kn:D1}, once $r$ is large enough, the image in projective space is \lq\lq balanced'' with respect to a suitable metric on the underlying vector space.
\item By the standard Kemp-Ness theory and the elucidation by Phong and Sturm \cite{kn:PS} (see also the related results of Luo \cite{kn:L} and Zhang \cite{kn:Z}) of the constructions of \cite{kn:D1}, the balanced condition implies Chow stability. 
\end{itemize}

 We have all the background in place to state the main result of this paper.
 Let $\cX$ and $\cX'_{p}$ be degenerations of $X$, as considered above, with a $\bC^{*}$-action on the central fibre $W'_{p}\subset \cX'_{p}$ generated by $A'_{p}$. Let $N(A'_{p})$ be the difference between the maximum and minimum eigenvalues of $A'_{p}$.

\begin{thm}
If $(X,L)$ is bi-asymptotically stable and $\cX$ satisfies the hypotheses above then there is a constant $K>0$ such that for infinitely many $p$ we have
$$   \Ch(W'_{p}) \geq K p^{-1} N(A'_{p}). $$
\end{thm}

The point of this result is that the $b$-stability of $(X,L)$ requires that an inequality of this kind holds,  and in fact the definition of $b$-stability is precisely that  similar inequalities hold for a more general class of  degenerations of $X$. Thus the theorem can be seen as a step towards the proof that the existence of a constant scalar curvature metric implies $b$-stability.
Said in another way, the theorem can be seen as a new obstruction to the existence of constant scalar curvature metrics. Suppose, say, we have an $\cX$ as above such that 
     $  p \ \Ch(W'_{p})/N(A'_{p})\rightarrow 0$ as $p\rightarrow \infty$. Then if the automorphism group of $(X,L)$ is finite we deduce that $X$ does not have a constant scalar curvature metric in the class $c_{1}(L)$.

\

We close this Introduction by repeating the algebro-geometric question raised in \cite{kn:D2}--regarding the relevance of the notion of b-stability. It might happen  that for some $p$ the component $R$ is contracted by the birational map and $W'_{p}$ is  irreducible. In that case, taking the higher power of the line bundle and replacing $W$ by $W'_{p}$, one can avoid discussing the birational modifications further and there should  be implications for the understanding of Gromov-Hausdorff limits of K\"ahler-Einstein metrics. So we ask 
\begin{question}
Is there an example of an equivariant degeneration $\cX$ as above, with $X$ Fano and $L=K_{X}^{-m}$ such that $R$ is not contracted in any $\cX'_{p}$?
Conversely, are there some hypotheses on $X$ which imply that $R$ always will be contracted, for some $p$?
\end{question} 
 
 See the further discussion in Section 3 below.
 
 \
 
 The author is grateful to Richard Thomas for many very helpful discussions.

\section{The set-up}

Consider an equivariant flat family $(\cL,\cX, \pi)$, with $\bC^{*}$-action, as in the previous section, and the inclusion maps $i_{p,\mu}$. 
\begin{lem}
If $i_{p,\mu}$ is not surjective then $\mu \leq C p$ where
$C= {\rm degree}(D)/{\rm degree}(R)$.
\end{lem}
By assumption there is a meromorphic section $s'$ of $\cL^{p}$ with a pole of order exactly $\mu$ along $R$. Thus $s=t^{\mu} s'$ is holomorphic, the restriction of $s$ to $R$ is not zero and $s$ vanishes to order $\mu$ on $B$. So the restriction of $s$ to $R$ is a non-zero section of $\cO(p)$ vanishing to order at least $\mu$ on $D$ and the assertion follows from basic facts about the degree. 

(Note that here, and in what follows, we use the standard convention that $t$ is the function on $\cX$ which is, logically speaking, the same as $\pi$.)

\

This lemma means that $m(p)$ is well-defined and in fact $m(p)\leq C p$. 
Now fix $m=m(p)$ and consider the sheaves $\pi_{*}(\cL^{p}), \pi_{*}(\cL^{p}\otimes \Lambda^{(m)})$ over $\bC$. They are torsion free and hence locally free.
\begin{lem}
There are sections $\sigma_{a}, \ (a=1,\dots ,N)$ which form a basis for $\pi_{*}(\cL^{p})$ and positive integers $\mu_{a}\leq m$ such that $\osigma_{a}= t^{-\mu_{a}} \sigma_{a}$ form a basis for $\pi_{*}(\cL^{p}) \otimes \Lambda^{(m)})$.
\end{lem}

This follows from basic facts about maps between modules over a PID, applied to the inclusion map $\pi_{*}(\cL^{p})\rightarrow \pi_{*}(\cL^{p}\otimes \Lambda^{(m)})$. The statement is also the content of Lemma 4 in \cite{kn:D2}. It follows from the definitions that $\sigma_{a}$ vanishes to order exactly $\mu_{a}$ along $B$ and $\osigma_{a}$ has a pole of order exactly $\mu_{a}$ on $R$. 
\begin{lem}
We can choose $\sigma_{a}$ to be eigenvectors for the induced $\bC^{*}$-action on $\pi_{*}(\cL^{p})$. 
\end{lem}
 Following the development in \cite{kn:D2}, we start with a flag in $H^{0}(W;\cL^{p}\vert_{W})$ defined by the maximal order of vanishing on $B$ of extensions over $\cX$. This flag is clearly $\bC^{*}$-invariant, so we can choose a compatible basis of eigenvectors in $H^{0}(W;\cL^{p}\vert_{W})$. Then any extension of one of these, vanishing to maximal order on $B$, gives a choice of $\sigma_{a}$. By a standard argument, averaging over the action of $S^{1}\subset \bC^{*}$, we can choose the extensions to be eigenvectors. 

Now write $\lambda_{a}$ for the weight of the  $\bC^{*}$-action on $\sigma_{a}$, so $g^{*}_{t}(\sigma_{a})= t^{\lambda_{a}} \sigma_{a}$.
 The projective embedding of $\cX$ defined by the sections of the line bundle $\cL^{p}$ can be written down as follows. Let $A_{p}$ be the diagonal matrix with diagonal entries $\lambda_{a}$ and let $X_{t}= t^{A_{p}} X$ where $X$ is embedded using the restrictions of the $\sigma_{a}$ to $\pi^{-1}(1)$. Then $\cX$ is embedded in $\bP^{N-1}\times \bC$ as the closure of
$$   \bigcup_{t\neq 0}(t,X_{t}). $$
This just expresses the definition, in terms of our basis $\sigma_{a}$. (It corresponds to the composite of the original embedding, defined by sections of $\cO(1)$, and the Veronese embedding of degree $p$.) In just the same way the image of the rational map defined by the sections of
$\cL^{p}\otimes \Lambda^{m}$ is given by replacing $X_{t}$ by $X'_{t}= t^{A'_{p}} X$ where $A'_{p}$ is the diagonal matrix with diagonal entries $\lambda_{a}-\mu_{a}$. This just expresses the definition, in terms of our basis $\osigma_{a}$. On the other hand it clearly reproduces the definition in \cite{kn:D2}.

Notice that throughout the construction we can allow the special case when
$m(p)=0$
 In that case we just take
$\cX'_{p}$ to be $\cX$, embedded by $\cL^{p}$.

To clarify notation, let $e_{a}$ be the standard basis of $\bC^{N}$ and $e^{*}_{a}$ the dual basis. Thus the $e^{*}_{a}$ are sections of $\cO(p)$ and the $e^{*}_{a}$ with $\mu_{a}>0$ form a basis for the sections which vanish on $B$. 

Let $\bP_{\infty}\subset \bP^{N}$ be the projectivization of the vector subspace spanned by vectors $e_{a}$ with $\mu_{a}>0$.

Recall that $W_{p}'$ is the central fibre of $\cX'_{p}$, so there is a rational map
$   j_{p}: W\rightarrow W'_{p}$.
By definition this is regular on the subset $\Omega=B\setminus D$ of $W$ and defines a $\bC^{*}$-equivariant open embedding of $\Omega$ in $\Omega'_{p}=W'_{p}\cap(\bP^{N}\setminus \bP_{\infty})$.
 Let  $B'_{p}$ be the closure in $W'_{p}$ of the image, so $B'_{p}$ is an irreducible component of $W'_{p}$. Let $R'_{p}$ be the union of all other irreducible components of $W'_{p}$.

We can study  $W'_{p}$ using \lq\lq arcs'' in $\cX$. Let $\Gamma:\Delta\rightarrow \cX$ be a holomorphic map from a disc with $\pi\Gamma(s)=0$ if and only if $s=0$.  Then $\Gamma$ induces a map $\tilde{\Gamma}$ from $\Delta$ to $\cX'_{p}$ and we can define an equivalence relation by $\Gamma_{1}\sim_{p}\Gamma_{2}$ if $\tilde{\Gamma_{1}}(0)=\tilde{\Gamma_{2}}(0)$, so points of $W'_{p}$ are equivalence classes of arcs. Fix a point $x\in W$ and a  local trivialisation of $\cL$ around $x$. Consider an arc $\Gamma$ with $\Gamma(0)=x$. With respect to this trivialisation the composites $\osigma_{a}\circ \Gamma$ are meromorphic functions on the disc. Then $\Gamma$ defines a point in $\Omega'_{p}=W'_{p}\cap(\bP^{N}\setminus \bP_{\infty})$ if, for each $a$, the function  $\osigma_{a}\circ \Gamma$ is bounded, and two such paths $\Gamma_{1}, \Gamma_{2}$ define the same point in $\Omega'_{p}$ if $\osigma_{a}\circ\Gamma_{1}(0)= \osigma_{a}\circ\Gamma_{2}(0)$ for all $a$. 
\begin{lem} Suppose that the sections of $\cO(p)$ over $\cX$ generate 
  $\cI_{B}(p)$, where $\cI_{B}$ is the ideal sheaf of $B$ in $\cX$ and that
  $x$ is a point in $B$. Let $\Gamma_{1}, \Gamma_{2}$ be  arcs through $x$ defining the same point  of $\Omega'_{p}$. Then for any holomorphic function $\tau$  defined on a neighbourhood of $x$  and vanishing to order
$\nu$ on $B$ the composites $t^{-\nu} \tau\circ \Gamma_{1}, t^{-\nu} \tau\circ\Gamma_{2}$ are bounded and have the same limit at $s=0$.
\end{lem}
First suppose that $\nu=1$. The hypothesis that $\cO(p)$ generates the ideal sheaf means that we can write $\tau=\sum_{\mu_{a}>0} f_{a} \sigma_{a}$ where $f_{a}$ are holomorphic. Thus 
$$t^{-1}\tau= \sum_{\mu_{a}>0} f_{a} (t^{\mu_{a}-1} \osigma_{a}), $$
and $$  t^{-1}\tau\circ \Gamma_{i}= \sum_{\mu_{a}>0} f_{a} (t^{\mu_{a}-1} \osigma_{a}\circ \Gamma_{i}),$$
from which the statement is clear. Now suppose that $\nu>1$. The global hypothesis
(6) on $B$ (from Section 1) asserts that $\tau$ is in $\cI_{B}^{\nu}$ and so can be written as a sum of terms which are products of the form $\tau_{1}\tau_{2}\dots \tau_{\nu}$ where $\tau_{j}$ vanishes on $B$. Then $t^{-\nu} \tau$ is a sum of products 
$$    (t^{-1}\tau_{1})(t^{-1}\tau_{2})\dots (t^{-1}\tau_{\nu}), $$
and the statement follows from the first case, applied to the $t^{-1} \tau_{i}$.

\begin{prop}
\begin{enumerate}
\item If the sections of $\cO(p)$ over $\cX$ generate $\cI_{B}(p)$ then
$R'_{p}$ is contained in $\bP_{\infty}$.
\item If the sections of $\cO(p)$ over $\cX$ generate $\cI_{B}(p)$ for $p=p_{1}$ and $ p=p_{2}$ then there is a $\bC^{*}$-equivariant isomorphism from $\Omega'_{p_{1}}$ to $\Omega'_{p_{2}}$, compatible with the open embeddings $j_{p_{i}}: \Omega \rightarrow \Omega'_{p_{i}}$.
\end{enumerate}
\end{prop}

(1) Suppose that $\Gamma$ is an arc defining a point in $W'_{p}$ and $x=\Gamma(0)$. If $x$ does not lie in $B$ there is (by the global generation hypothesis) a section $\sigma_{a}$ with $\lambda_{a}>0$ such that $\sigma_{a}$ does not vanish at $x$. Then $\osigma_{a}\circ \Gamma$ is unbounded and $\Gamma$ defines a point in $\bP_{\infty}$. So suppose that there is a  component of $R'_{p}$ which does not lie in $\bP_{\infty}$. Choose a point $z$ in $R'_{p}\cap(\bP^{N}\setminus \bP_{\infty}$  which does not lie in the closure of $B'_{p}$ and an  arc $\Gamma$ with $\Gamma(0)\in B$ representing $z$. There is a polynomial
$F$ which vanishes on $B'_{p}$ but not at $z$. Then $\tau= F(\osigma_{a} )$ is a meromorphic section which vanishes to some order $\nu>0$ on $B$.
The fact that $F$ does not vanish at $z$ means that $t^{-\nu} \tau\circ \Gamma$ is unbounded, contrary to  Lemma 4.

\

(2) Suppose $\Gamma_{1}$ and $\Gamma_{2}$ define the same point in $\Omega'_{p_{1}}$.
By considering only the holomorphic sections $\osigma_{a}$ (i.e. those with $\mu_{a}=0$), we see that $\Gamma_{1}(0)=\Gamma_{2}(0)$. Working in a local trivialisation around this point the second statement of  Lemma 4 shows that for any local section $\tau$ of a line bundle  vanishing to order $\nu$ on $B$ the composites $t^{-\nu} \tau\circ \Gamma_{i}$ are bounded and have the same limit. In particular this is true for the sections of $\cL^{p_{2}}$. So we see that $\Gamma_{1}\sim_{p_{1}}\Gamma_{2}$ implies that $\Gamma_{1}\sim_{p_{2}} \Gamma_{2}$. Reversing the roles of $p_{1}, p_{2}$ we see the converse and we  get a bijection between $\Omega'_{p_{1}}$ and $\Omega'_{p_{2}}$. It is clear from the construction  that this is an equivariant isomorphism.

 From now on we consider sufficiently  large values of $p$, say $p\geq p_{0}$ so that the global generation hypothesis is satisfied.

Now we turn to consider the $\bC^{*}$-actions on $\cX$ and $\cX'_{p}$.
There is no loss in supposing that the weights of the action on sections of $\cO(1)$ over $W$ run from $0$ to $-M$ say, where $M>0$. Then the weights $\lambda_{a}$  of the action on sections of $\cO(p)$ over $W$ run from $0$ to $-Mp$. 
\begin{lem}
We can choose a basis $\sigma_{a}$ as above and an equivariant lift  $\iota:\bC\rightarrow \cX$ of $\pi$ with image $\Sigma\subset \cX$ such that
\begin{itemize}
\item  $x_{0}=\iota(0)$ is in $ B$; 
 \item $\lambda_{0}=0$;
\item $\sigma_{0}$ does not vanish anywhere on $\Sigma$ but for $a>0$ all $\sigma_{a}$ vanish on $\Sigma$. 
\end{itemize}
\end{lem} 

To see this, consider the restriction of the $\bC^{*}$ action to $B$. By assumption (5) $B$ does not lie in a linear subspace and this means that there must be points  of $B$ represented by vectors of highest and lowest weight. The point relevant to us is a point of highest weight. We can obviously choose the basis so that this is $x_{0}=[e^{*}_{0}]$.. In dynamical language, this is a repulsive fixed point for the $\bC^{*}$-action on the central fibre and there is a unique point $x_{1}$ in $X=\pi^{-1}(1)$ which flows to $x$ in the limit as $t\rightarrow 0$, i.e. $\lim_{t\rightarrow 0} g_{t}(x_{1})=x_{0}$. The $\bC^{*}$-orbit of $x_{1}$ defines the lift $\iota$.

\

Since $\sigma_{0}$ does not vanish at $s(0)\in B$ me must have $\mu_{0}=0$.
For any $p>0$ the lift $\iota$ induces a lift $\iota':\bC\rightarrow \cX'_{p}$ with image $\Sigma_{p}$ and $\iota'(0)$ is again a highest weight vector for the action on $\cX'_{p}$. Since all the sections $\sigma_{a}$ with $\mu_{a}>0$ vanish on $\Sigma$  this $\iota'(0)$ certainly does not lie in $\bP_{\infty}$.  Thus (at least when $p$ is sufficiently large) the point $\iota'(0)$ lies in $\Omega_{p}\subset B'_{p}$. It is clear that if we have two values $p_{1}, p_{2}$ the isomorphism  from $\Omega'_{p_{1}}$ to $\Omega'_{p_{2}}$ is compatible with these points and extends  to an isomorphism between neighbourhoods $U_{p_{i}}$  of the $\Sigma_{p_{i}}$ in $\cX'_{p_{i}}$.

\section{Examples and discussion}

  Work in affine space $\bC^{3}$ and consider the $\bC^{*}$-action
$$    (x,y,z)\mapsto (t^{-\alpha} x, t^{-\beta} y , t^{-\gamma}z)$$
If $\alpha,\beta,\gamma>0$ then the origin is a repulsive fixed point of the kind we are considering. Let $Z$ be the affine variety given by the equation $ xy= z$. Then the transform $Z_{t}$ under the action is given by the equation $xy = t^{\gamma-\alpha-\beta} z$.  If $\gamma>\alpha+\beta$ the limiting variety
$Z_{0}$ is the union of two planes. Suppose for definitenes that $\gamma-\alpha-\beta=1$.Let $B$ be the plane $y=0$, say. The total space ${\cal Z}\subset \bC^{4}$ is the quadric cone with equation $xy=tz$. The lift $\Sigma$ is just given by $x=y=z=0$. The plane $B\subset {\cal Z}$ defined by $y=t=0$ is not a Cartier divisor but one can check that it does satisfy our technical condition (6). This example does not quite fit our framework because $B$ does lie in a proper linear subspace. But if take $\alpha=2, \beta=3, \gamma=10$ and start instead with the smooth surface $Z$ in $\bC^{3}$ given by the equation $(x^{3}-y^{2}) y= z$ the total space $\cZ$ has equation $(x^{3}-y^{2}) y = t z$, the central fibre $W$ is the union of $B$, defined by $x^{2}-y^{3}=0$, and $R$, defined by $y=0$. Again one can check that $B$ satisfies our technical condition (6). Take $p=3$: the space   of cubics vanishing on $B$ is spanned by $x^{2}-y^{3}$. To construct $\cZ'_{3}$ we can embed $\cZ$ in $\bC^{5}$ with co-ordinates
$(x,y,z,u,t)$ defined by the equations $u=x^{3}-y^{2}, uy=tz$. Now we replace $u$ by $ut$ to get equations $ut=x^{3}-y^{2}, u y=z$ defining $\cZ'_{3}\subset \bC^{5}$. Then $B'_{3}\subset \bC^{4}$ is defined by the equations $x^{3}-y^{2}=0, uy=z$. The line $x=y=z=0$ lies in $B'_{3}$ and and is contracted under the birational equivalence with $B$.

\

We can now raise a question which can be viewed as a  local version of Question 1. As discussed briefly in \cite{kn:D2}, Question 1 is related to the following.
Let $\Lambda^{*}$ denote the sheaf of functions on $\cX$ with arbitrary poles on $R$. Then the sections of $\cL^{p}\otimes \Lambda^{*}$, for all $p\geq 0$, define a graded algebra over the ring of holomorphic functions over $\bC$ and the issue is whether this is finitely generated. A local version of this is to ask whether the sheaf $\Lambda^{*}$ is finitely generated as a sheaf of algebras over $\cO_{\cX}$. Let $\cI^{(\mu)}_{B}$ denote the sheaf of holomorphic functions on $\cX$ which vanish to order $\mu$ on $B$ (the \lq\lq symbolic power'' of $\cI_{B}$). Then we can form
$$  \cI_{B}^{*}= \bigoplus_{\mu} \cI^{(\mu)}_{B}. $$
This is a sheaf of algebras and there is a surjective sheaf homomomrphism
$$  T: \cI_{B}^{*}\rightarrow \Lambda^{*}$$ defined by $ T(f)= t^{-\mu} f$ for $f\in \cI_{B}^{(\mu)}$. Now examining the proof of Lemma 4 above we see that the essential thing we need is that $\Lambda^{*}$ is finitely generated, as a sheaf of algebras over $\cO_{\cX}$. Since $T$ is surjective this will be true if $\cI_{B}^{*}$ is finitely generated and in turn this will certainly be true if $\cI_{B}^{(\mu)}$ coincides with the power $\cI_{B}^{\mu}$, which we assumed in our situation.

 Examining the proofs further we see that we only need this finite generation property at a special point of $B$. So in sum we arrive at the following
\begin{question}
Let $Z\subset \bC^{N}$ be a smooth affine variety containing $0$. Let
$g_{t}:\bC^{N}\rightarrow \bC^{N}$ be a linear $\bC^{*}$-action such that $0$ is a repulsive fixed point. For $t\neq 0$ let $Z_{t}=g_{t}(Z)$ and let $\cZ\subset \bC^{N}\times \bC $ be defined as the closure of the family. Suppose that $B$ is a reduced component of $Z_{0}=\cZ\cap (\bC^{N}\times \{0\})$ and $R$ is the union of all other components. Let $\Lambda^{*}_{0}$ be the germs of meromorphic functions on $\cZ$ around $(0,0)\in \bC^{N}\times \bC$ with poles along $R$. Then is $\Lambda^{*}_{0}$ finitely generated as an algebra over the germs of holomorphic functions?
\end{question}

If the answer to this question is positive it would be possible to omit our technical hypothesis (5) in proving Theorem 2. On the other hand if the answer is negative it seems there could be an essential difficulty in extending this approach.

 \section{Blowing up}

We follow the approach of Stoppa \cite{kn:S1}, \cite{kn:S2}. We begin by reducing to a simpler situation. Starting with $\cX$ and $p$ we define
$\cX'_{p}$ and a divisor $R'_{p}\subset \cX'_{p}$. Now given $q$ we can make the same construction to define a rational map from $\cX'_{p}$ to the product of a projective space with $\Delta$ using meromorphic sections of $\cO(q)$ over $\cX'_{p}$ with poles along $R'_{p}$. Denote the result by $\left( \cX'_{p}\right)'_{q}$.
Then we have a canonical isomorphism
$$   \left(\cX'_{p}\right)'_{q}= \cX'_{pq}. $$
This was stated, with outline proof, in \cite{kn:D2}. We can give a proof more  in   line with our current point of view as follows. A meromorphic section of $\cO(q)$ with poles along $R'_{p}$ can be written as $t^{-\nu} \tau$ where $\tau$ is a holomorphic section of $\cO(q)$ vanishing to order $\nu$ on $B'_{p}$.
The hypothesis that $H^{0}(X,L)$ generates $\bigoplus H^{0}(X,L^{k})$ means that we can write $\tau= \overline{P}(\osigma_{a})$ where $\overline{P}$ is a homogeneous polynomial of degree $q$ with coefficients meromorphic in $t$. Thus with a slight stretch of notation we can also write
$\tau=P(\sigma_{a})$ where $P$ is another polynomial with meromorphic coefficients. That is to say we simply substitute $\osigma_{a}= t^{-\mu_{a}} \sigma_{a}$. Thus $\tau$ can  be regarded as a meromorphic section of $\cO(pq)$ over $\cX$ vanishing to order $\nu$ on $B$. Thus $t^{-\nu} \tau$ is also one of the sections we use in defining $\cX'_{pq}$. Conversely if we start from a section
$\tau$ of $\cO(pq)$ over $\cX$ vanishing to order $\nu$ on $B$ the fact that the sections of $L^{p}$ over $X$ generate those of $L^{pq}$ implies that we can write $\tau$ as a polynomial in the $\sigma_{a}$ with meromorphic coefficents. This gives rise to the equivalence of $\cX'_{pq}$ and $\left(\cX'_{p}\right)'_{q}$.

Recall that we can fix $p_{0}$ so that for all $p\geq p_{0}$ the sections $\Sigma_{p}\subset \cX'_{p}$ have isomorphic neighbourhoods. If we restrict attention to multiples $p=qp_{0}$ we can think of $\cX'_{p}$ as being obtained from $\cX'_{p_{0}}$ by meromorphic sections of $\cO(q)$, as above. Then the section $\Sigma'_{p_{0}}$ does not meet $R'_{p_{0}}$.  With this said, there is no real loss of generality in supposing that in fact $\Sigma\subset \cX$ {\it does not meet} $R$, i.e. that
$\iota(0)$ is in $\Omega$. In this case there are, for all $p$, neighbourhoods $U_{p}$ of the $\Sigma_{p}$ in $\cX'_{p}$ which are isomorphic to a neighbourhood $U$ of $\Sigma$ in $\cX$. We will make this assumption from now on to simplify  the notation---but it is no restriction because we can always replace $\cX$ by $\cX'_{p_{0}}$ in what follows.

    We blow up $\Sigma$ in $\cX$ to get a variety $\hcX$ which is yields another equivariant flat family $\hpi:\hcX\rightarrow \bC$. The ideal sheaf $\cI_{\Sigma}$ corresponds to a line bundle $\cE\rightarrow \hcX$ and for suitable values of $r,\gamma$ the line bundle $\cL^{r\gamma}\otimes \cE^{r}$ will define a projective embedding of $\hcX$. The non-zero fibres of $\hpi$ are just given by blowing up $g_{t}(X)$ at the point $\iota(t) $ but, as Stoppa explains, the fibre $\hpi^{-1}(0)$ is not necessarily the blow up of $W$ at $x_{0}$. In general it may contain another component, or union of components $P$.

\

To illustrate this phenomenon, consider the first example in Section 3, so we want to blow-up the quadric cone $\{ xy=tz\}$ along the line $x=y=z=0$.
The blow up $\hat{\cZ}$ is a subvariety of $\bC^{4}\times \bP^{2}$.  Let $\hat{Z_{t}}$ be the blow-up of $Z_{t}\subset \bC^{3}$ at the origin $x=y=z=0$.  Then $\hat{\cZ}$ is the closure of the union over non-zero $t$ of $(t,\hat{Z_{t}}$. Let $(r_{1}, r_{2}, r_{3})$ be a vector with $r_{1}r_{2} \neq 0$. The line in $\bC^{3}$ generated by this vector meets $Z_{t}$ at the point
$$   \frac{tr_{3}}{r_{1}r_{2}}( r_{1}, r_{2}, r_{3}) $$
As $t$ tends to zero this point tends to the origin. It follows that $\hat{\cZ}$ contains the whole of $(0,\bP^{2})\subset \bC^{4}\times \bP^{2}$. Thus $P$ is a copy of $\bC\bP^{2}$. 

\

In any case we write $\widehat{W'}$ for the central fibre of $\widehat{\cX}$.

Now consider $\cX'_{p}$ and the section $\Sigma_{p}\subset \cX'_{p}$. We can blow up $\Sigma_{p}$ to obtain $\widehat{\cX'_{p}}$.  Since blowing up is a local operation there are neighbourhoods $\widehat{U_{p}}\subset \widehat{\cX'_{p}}$ of the exceptional set which are isomorphic, for all $p$.

We want to consider projective embeddings of $\widehat{\cX'_{p}}$.  Recall that we  write $\pi^{(p)}: \cX'_{p}\rightarrow \bC$.

\begin{prop} Suppose $p,r$ are chosen so that the restriction map
    $$   \pi_{*}\cL^{p} \rightarrow \pi_{*}( \cL^{p}\otimes \cO_{r\Sigma}) $$
is surjective and the sections of $\cL^{p}\otimes \cI_{\Sigma}^{r}$ define a projective embedding of $\widehat{\cX}$. Then
the restriction map $\pi^{(p)}_{*}(\cL^{p} \otimes \Lambda^{(m(p))})\rightarrow \pi^{(p)}_{*}( \cL^{p}\otimes  \Lambda^{(m(p))}\otimes \cO_{r\Sigma_{p}})$ is surjective and the sections of $\cL^{p}\otimes \Lambda^{(m(p)}\otimes \cI_{\Sigma_{p}}^{r}$ define a projective embedding of $\widehat{\cX'_{p}}$.
\end{prop}

(Here of course we write $\cO_{r\Sigma}= \cO/\cI^{r}_{\Sigma}$.)

The surjectivity of the restriction map follows simply from the fact that
the sheaf $\Lambda^{(m)}$ contains the holomorphic functions. Now starting with the family $\widehat{\cX}$ we can apply the same construction using the functions with poles along the copy of $R$ in the central fibre of 
$\widehat{\cX}$.  This gives a projective family $\left(\widehat{\cX}\right)'_{p}$ say. But it is clear that this is canonically isomorphic to $\widehat{\cX'_{p}}$. In terms of $\cX$, we are in both cases considering the image of the rational map defined by sections of $\cL^{p}\otimes \cI_{\Sigma}^{r}\otimes \Lambda^{(m)}$.

All this can be summarised by saying that, under the assumption that $x_{0}$ is not in $R$, the two constructions (passing from $\cX$ to $\cX'_{p}$ and blowing up) do not interact. 

\section{Chow invariant calculations}

Suppose that $p,r$ are chosen as in Proposition 2. We have two projective varieties with $\bC^{*}$-action: the central fibre $W'$ of $\cX'_{p}$ and the central fibre $\widehat{W'}$ of $\widehat{\cX'_{p}}$. We want to compare the Chow invariants of the two. The main task is to compare the terms from the traces of the actions. This is essentially contained in Stoppa's work but we take a slightly different approach, which  admits extensions to more general (non equivariant) degenerations.  Let $\cE\rightarrow \widehat{\cX'_{p}}$ be the line bundle defined by the blow up and $Y\subset \widehat{\cX'_{p}}$ be the exceptional divisor. For any $r,k$ we have a line bundle $\cO(k)\otimes \cE^{r}$ over $\widehat{\cX'_{p}}$ and the direct image is an equivariant vector bundle over $\bC$. This is the same as taking the direct image of the sheaf $\pi^{(p)}_{*}(\cO(k)\otimes \cI_{\Sigma'_{p}}^{r})$, working over  $\cX'_{p}$. We write $\hTr(p,k,r)$ for the trace of the action on the central fibre of this bundle. Similarly we write
$\Tr(p,k)$ for the corresponding trace formed from the sections of $\cO(k)$ over $\cX'_{p}$. 
\begin{prop} There is function $F(r)$ such that if 
$\pi^{(p)}_{*}(\cX'_{p};\cO(k))$ maps onto $\pi^{(p)}_{*}(\cX'_{p}; \cO(k)\otimes \cO_{r\Sigma_{p}})$ then $\hTr(p,k,r)= \Tr(p,k)+ F(r)$. Furthermore, for large enough $r$ the function $F(r)$ is a polynomial of degree $n+1$ in $r$.  
\end{prop}

Suppose in general that $E$ is an $\bC^{*}$-equivariant vector bundle over $\bC$. The $\bC^{*}$-action defines a trivialisation of  $E$ away from $0$ and so the first Chern class of $E$ is defined, relative to this trivialisation. This integer coincides with the trace of the action on the central fibre. From this of view we have to compare the  Chern classes of the bundles $\pi^{(p)}_{*}(\cO(k))$ and $\pi^{(p)}_{*}(\cO(k)\otimes \cI_{\Sigma_{p}}^{r})$. The surjectivity  hypothesis of the proposition means that there is an exact sequence of bundles:
$$  0\rightarrow  \pi^{(p)}_{*}(\cO(k)\otimes \cI_{\Sigma_{p}}^{r} \rightarrow \pi^{(p)}_{*}(\cO(k))\rightarrow \pi^{(p)}_{*}(\cX'_{p}; \cO(k)\otimes \cO_{r\Sigma_{p}})\rightarrow 0, $$
so the difference of the Chern clases is the Chern class of
 $\pi^{(p)}_{*}(\cO(k)\otimes \cO_{r\Sigma_{p}}$. The fact that $\lambda_{0}=\mu_{0}=0$ means that the line bundle $\cO(1)$ is trivialised, as an equivariant bundle, over $\Sigma_{p}$ so we can omit the factor $\cO(k)$ and we see that
$\hTr(p,k,r)-\Tr(p,k)= c_{1}(\pi^{(p)}_{*}(\cO_{r\Sigma_{p}}))$ and this is plainly independent of $p$ since the $\Sigma_{p}$ have isomorphic neighbourhoods. So we have
$\hTr(p,k,r)-\Tr(p,k)=F(r)$.

The fact that $F$ is a polynomial of degree $n+1$ for large $r$ follows from Riemann-Roch for families. When $r$ is large enough there is an exact sequence
$$0\rightarrow \pi^{p)}_{*}(\cO/\cI_{\Sigma_{p}}^{r-1})\rightarrow \pi^{(p)}_{*}(\cO/\cI_{\Sigma_{p}}^{r})\rightarrow \varpi_{*}( \cE\vert_{Y} ^{r}) \rightarrow 0, $$
where $\varpi:Y\rightarrow \bC$ is the projection map from the exceptional divisor. Thus 
$F(r)-F(r-1)$ is the first Chern class of $\varpi_{*}( \cE\vert_{Y}^{r})$. Use the trivialisation defined by $g_{t}$ of $\cX$ over $\bC\setminus\{0\}$ to compactify $\cX$ to a family over $\bC\bP^{1}$. This defines a compactification of the exceptional divisor which is a compact variety of dimension $n$. When $r$ is large enough for the higher cohomology to vanish, the first Chern class is given by the Grothendieck-Riemann-Roch formula  as a polynomial of degree $n$ in $r$. Now summing over $r$ we see that $F$ is a polynomial of degree $n+1$, for large $r$.

Now we can compare $I(W'_{p})$ and $I(\hat{W'_{p}})$. Write
$$  F(r)= f r^{n+1} + \epsilon_{1}(r) $$
where $\epsilon_{1}(r)=O(r^{n})$. Fix $p,r$ and consider the projective embedding of $\widehat{\cX'_{p}}$ defined by $\cO(k)\otimes \cI_{\Sigma_{p}}^{kr}$. We see that
$$   I(\hat{W'_{p}})= \lim_{k\rightarrow \infty} k^{-n-1}\hTr(p,k,kr), $$
while $$ I(W'_{p})=\lim_{k\rightarrow \infty}k^{-n-1} \Tr(p,k), $$
so \begin{equation}I(\hat{W'_{p}})- I(W'_{p})= \lim_{k\rightarrow \infty} k^{-n-1} F(kr)= f r^{n+1}. \end{equation}

Now fix $p,r$ such that sections of $\cO(p)$ over $\cX$ generate $\cO_{r\Sigma}$.
Write $\Tr=\Tr(p,1)$  and $\hTr= \hTr((p,1,r)$. These are the trace terms corresponding to the actions on $W'_{p}$ and $\widehat{W'_{p}}$ respectively.
Similarly write $\Dim, \hDim$ for the dimensions of the underlying spaces in the projective embeddings of $W'_{p}, \widehat{W'_{p}}$. By considering the restriction to a general fibre we see that
$\Dim-\hDim$ is the dimension of the space of poloynomials of degree at most $r$ in $n$ variables, i.e.
$$  \Dim-\hDim= \frac{ (r+1)(r+2)  \dots (r+n)}{n!}$$
Thus $$\hDim= \Dim-(d_{0}r^{n}+ d_{1}r^{n-1}+ \epsilon_{2}(r))
$$
where $d_{0}, d_{1}>0$ and $\epsilon_{2}(r)=O(r^{n-2})$.

Let $\Vol, \hVol$ be the degrees (or volumes) of $W_{p}', \widehat{W_{p}'}$ and let $I=I(W'_{p}) ,\hat{I}= I(\hat{W'_{p}})$.  Then $\Vol,\hVol$ are given by the degrees of $X$ and of the blow-up of $X$ at a point, so we get
$  \Vol= V p^{n}$ for some $V>0$ and
$\hVol= \Vol-d_{0}r^{n}$. Finally we know that $\Dim$ is the dimension of $H^{0}(X,L^{p})$ so
$$  \Dim=V p^{n}+ O(p^{n-1}). $$
Putting everything together gives the following
\begin{prop} We have  \begin{equation} \Ch(W_{p}')= \frac{\Tr}{\Dim}- \frac{I}{V p^{n}},\end{equation} 
    \begin{equation} \Ch(\wW'_{p})=\frac{ \Tr- f r^{n}+\epsilon_{1}(r)}{\Dim-(d_{0}r^{n}+ d_{1} r^{n-1}+\epsilon_{2}(r)}- \frac{I- f r^{n+1}}{V p^{n}- d_{0}r^{n}} \end{equation}
where $\epsilon_{1}(r)$ is $O(r)$, $\epsilon_{2}(r)$ is $O(r^{n-2})$, $\Dim =Vp^{n}+ \eta(p)$ with $\eta(p)=O(p^{n-1})$ and $ d_{1}>0$.

\end{prop}

Now let $p=r\gamma$ for an integer $\gamma$.

\begin{cor}
We can write
$$\Ch(\wW'_{r\gamma})= (1-d_{0}V^{-1}\gamma^{-n}) \Ch(W') + \frac{\Tr}{\Dim}\left( d_{1}r^{-1}\gamma^{-n}+\alpha(r,\gamma)\right) + \beta(r,\gamma), $$
where $\vert \alpha(r,\gamma)\vert \leq c (r^{-2}\gamma^{-n}+ r^{-1} \gamma^{-n-1})$ and $\vert \beta(r,\gamma)\vert\leq c (r^{-1}\gamma^{1-n}+ \gamma^{-n})$.
\end{cor}
This just uses the information contained in the preceding Proposition. We first observe that the term $\epsilon_{1}$ in (5) can be absorbed in $\beta$. Similarly, the term $\epsilon_{2}$ can be absorbed in $\alpha$, using the fact that $\Dim\sim V(r\gamma)^{n}$. Next one sees that the two terms involving $f$ in (5) cancel, modulo terms that can be absorbed in $\beta$. Thus we are reduced to considering the simpler expression
$$  Q= \frac{\Tr}{\Dim-(d_{0} r^{n}+d_{1} r^{n-1})} - \frac{I}{V r^{n} \gamma^{n}- d_{0}r^{n}}. $$
Write this as
$$  Q= \left( \frac{\Tr}{\Dim}- \frac{I}{V r^{n} \gamma^{n}}\right) (1-d_{0}V^{-1}\gamma^{-n})
+\frac{\Tr}{\Dim}\left( \frac{1}{1- (d_{0}r^{n}+d_{1}r^{n-1})/\Dim} - \frac{1}{1-d_{0}r^{n}/(V r^{n} \gamma^{n})}\right). $$
Then some elementary manipulation  and estimation gives the result.

We can now prove Theorem 2. Suppose that $X$ is $(\gamma,r)$-stable. This means that the Chow invariant of $\wW'$ is positive so we get
$$  \Ch(W') \geq (1-d_{0}\gamma^{-n}) \left( \frac{-\Tr}{\Dim}(d_{1}r^{-1} \gamma^{-n}+\alpha) + \beta\right). $$
Now recall that $\Tr$ is the trace of the action corresponding to $\cX'_{r\gamma}$ and this is {\it  less than} the corresponding trace for the action on $H^{0}(W,\cL\vert_{W}^{r\gamma})$ since we modify the action by subtracting the positive terms $\mu_{a}$. For the fixed variety $W$ we know by the standard asymptotics that the trace is given by a Hilbert polynomial of degree $n+1$. So we see that
$$  \frac{\Tr}{\dim}\leq - w(r\gamma), $$
for some $w>0$. Choose  $\gamma_{0},r_{0}$ so that if $\gamma\geq \gamma_{0}$ and $r\geq r_{0}$ we have \begin{itemize}
\item  $(1-d_{0}V^{-1} \gamma^{-n})\geq 1/2$ \item  $\alpha(r,\gamma)\leq \frac{d_{1}}{2} r^{-1} \gamma^{-n}$
\item $\beta(r,\gamma)\leq \frac{w d_{1}}{8}\gamma^{1-n}$ 
\end{itemize}
Then we get
$$\Ch(W') \geq \frac{1}{4} ( w d_{1}\gamma^{1-n}) - 2 \beta)\geq \frac{3}{16} w d_{1} \gamma^{1-n}. $$
Next recall  that we defined $N(A_{r\gamma}')$ to be the difference between the maximum and minumum eigenvalues of the generator  $A'_{r\gamma}$ of the action. This is
$$   N(A_{r\gamma}')= \max_{a} (-\lambda_{a}+\mu_{a})$$ and by Lemma 1 we see that
$N(A_{r\gamma}')\leq (C+M)r\gamma$. So we have 
$$ \Ch(W'_{r\gamma}) \geq K (r\gamma)^{-1} N(A_{r\gamma}')\gamma^{1-n},$$
where $K=\frac{3}{16}  (C+M) w d_{1}$. 
To sum up: suppose $(X,L)$ is bi-asymptotically stable. We first fix a large $\gamma$ as above and such that $(X,L)$ is $(\gamma,r)$-stable for sufficiently large $r$. Then as $p$ runs over sufficiently large multiples of $\gamma$ we have obtained the inequality stated in Theorem 1.



\end{document}